\documentclass[a4paper,12pt]{article}

\usepackage{amsmath,amsfonts}
\usepackage{mathrsfs} 
\usepackage{parskip}
\usepackage{enumerate}
\usepackage{graphicx}
\usepackage{bbm}
\usepackage{bm}

\newtheorem{theorem}{Theorem}

\def \Liminf{\mathop{\underline{\lim}}\limits}

\def\BB{\mathbb{B}}
\def\CC{\mathbb{C}}
\def\DD{\mathbb{D}}
\def\Pb{\mathbf{P}}

\def\Ex{\mathbf{E}}

\def\TT{\mathbb{T}}

\def\KK{\mathbb{K}}

\def\II{\mathbb{I}}
\def\JJ{\mathbb{J}}

\def\1{\mbox{1\hspace{-.25em}I}}
\begin{document}
\title{On  Multi-step MLE-process for Markov  Sequences}

\author{Yu. A. \textsc{Kutoyants}\\
{\small Universit\'e du Maine, Le Mans, France and}\\
%
%
{\small National Research University ``MPEI'', Moscow, Russia}\\[8pt]
A. \textsc{ Motrunich}\\
{\small Universit\'e du Maine, Le Mans, France and}\\
{\small Stat\'esia, Le Mans, France}
}

\date{}
\maketitle
\

\begin{abstract}
We consider the problem of the construction of the estimator-process of the
unknown finite-dimensional parameter in the case of the observations of
nonlinear autoregressive process. The estimation is done in two or three
steps. First we estimate the unknown parameter by a learning relatively short
part of observations and then we use the one-step MLE idea to construct
an-estimator process which is asymptotically equivalent to the MLE. To have
the learning interval shorter we introduce the two-step procedure which leads
to the asymptotically efficient estimator-process too. The presented results
are illustrated with the help of two  numerical examples. 

\end{abstract}

\section{Introduction}

This work is devoted to the problem of finite-dimensional parameter
estimation in the case of observations of  Markov sequence in the asymptotics
of large samples. The observations are  $X^n=(X_0,X_1,X_2,\ldots,X_n)$. For
simplicity of exposition we take as a
model of observations a nonlinear time series satisfying the relation 
\begin{equation}
\label{0}
X_j=S\left(\vartheta ,X_{j-1}\right)+\varepsilon _j,\quad j=1,2,\ldots
\end{equation}
and the initial value $X_0$ is given too. The random variables
$\left(\varepsilon _j\right)_{j\geq 1}$ are i.i.d. with some known smooth
density function $g\left(x\right)$. The function $S\left(\vartheta ,x\right)$
is supposed to be known and smooth with respect to $\vartheta $.  It can be
verifies that under the supposed regularity conditions the family of measures
corresponding to these model of observations is locally asymptotically normal
(LAN).  Our goal is to construct a sequence (we say process) of estimators
$\vartheta _n^\star=\left(\vartheta _{k,n},k=N+1,\ldots,n\right)$, where $N
\ll n$. By the first $N+1$ observations $X^N=\left(X_0,X_1,\ldots,X_N\right)$
we estimate the parameter $\vartheta $ and the obtained {\it preliminary}
estimator $\bar\vartheta _N$ we use in the construction of the estimator
process $\vartheta _n^\star$. This construction is based on the modification
of the well-known one-step maximum likelihood estimator (MLE) procedure
introduced by Le Cam in 1956 \cite{LC56} for LAN families of distributions.
In the proofs we follow the similar work \cite{Kut15} devoted to parameter
estimation in the case of ergodic diffusion process. Such estimator-processes
appeared in the works devoted to the problem of approximation of the solution
of backword stochastic differential equations (see review in \cite{K14}).  As
the initial estimator is constructed by a relatively small number of
observations $N \sim n^\delta $ with $\delta <1$ the rate of convergence of
the preliminary estimator is ``bad'' $\sqrt{N}\sim n^{\delta /2} $
$$
\sqrt{N}\left(\bar\vartheta _N-\vartheta
\right)\Rightarrow {\cal N}\left(0,\DD\left(\vartheta \right)\right)
$$ 
and we have {\it to improve} this rate up to the {\it optimal} $\sqrt{n}$
and {\it to improve} the limit variance up to the {\it optimal}.

Therefore this work is devoted to adaptive estimation for LAN family
\cite{FH82}. The structure of our estimator-processes is  in some sense close
to that of the 
Fisher-scoring algorithm, but the proposed realization is different because we
have to improve the rate of convergence. 
   The idea to use a preliminary estimator with a ``bad'' rate of
convergence in the one-step MLE framework to obtain asymptotically efficient
one was used by Skorohod and Khasminskii \cite{SKh96} and the idea to improve
the rate of convergence of preliminary estimator using multi-step
Newton-Raphson procedure was realized by Kamatani and Uchida \cite{KU15}. In
the work \cite{SKh96} it was considered the problem of parameter estimation
for partially observed diffusion processes and in \cite{KU15} it was
considered the problem of parameter estimation by the discrete time
observations of the diffusion process in the asymptotics of {\it high
  frequency observations} , i.e., they supposed that the step of
discretization tends to zero. We consider the multi-step pprocedure of
one-stem MLE type for Markov sequences. Another particularity of the presented work is
the following. We propose a sequence of estimators, which can be easily
calculated and the same time it has the same asymptotic properties as the
asymptotically efficient MLE. This means that these estimators are
asymptotically normal and that its limit variance is the inverse Fisher
information matrix.

The properties of the parameter estimators for nonlinear time series and
Markov sequences, of course, are well-known. Let us mention here the works by
Roussas \cite{RG65}, Ogata and Inagaki \cite{OI77}, Varakin and Veretennikov
\cite{VV02}). More about statistical problems for time series can be found in
the monographs by Veretennikov \cite{Ver00}, Taniguchi and Kakizawa
\cite{TK00}, Fan and Yao \cite{FY03}, and the references therein.

Note that we take the time series \eqref{0} just for simplicity of
expositions. The proposed results can be generalized on the more general
Markov sequences defined by their transition density if we suppose that this
density satisfies to the corresponding regularity conditions.

The process $\left(X_j\right)_{j\geq 0} $ has a transition density
\begin{align*}
\pi \left(\vartheta ,x,x'\right)=g\left(x'-S\left(\vartheta ,x\right)\right).
\end{align*}
It  depends on the
parameter $\theta$ and defines the probability of reaching the state $x'$
after sojourning in the state $x$. The parameter $\vartheta $  takes its
values in some open, convex,  bounded set $\Theta\subset 
R^d$.

The construction of the one-step MLE-process in this work is done in two
steps. On the first step we estimate the unknown parameter by the observations
$X^N=\left(X_0,X_1,\ldots,X_N\right)$ 
on the {\it learning interval} $j\in \left[0,N\right]$. As preliminary
estimator we can take the MLE, Bayes estimator (BE), estimator of the method
of moments (EMM) or any other estimator, which is consistent and
asymptotically normal.

Let us recall some of them. The MLE  is defined as follows. 
Introduce the  likelihood function 
\begin{equation}\label{1}
V(\vartheta ,X^n)=\pi (\vartheta ,X_0)\prod_{j=1}^n \pi(\vartheta
,X_{j-1},X_j), \qquad \vartheta\in \Theta .
\end{equation}
 We  suppose that the observations are strictly stationary and therefore  the
 density of the initial value is the density of the invariant measure
$\pi\left(\vartheta ,x\right)$.

The maximum likelihood estimator we introduce as usual by the equation 
\begin{equation}
\label{2}
V(\hat{\vartheta} _n,X^n)=\sup_{\vartheta \in\Theta}V(\vartheta,X^n).
\end{equation}
If this equation has many solutions then we can take any of them as the MLE.

It is known that under the regularity conditions the MLE  is consistent and asymptotically normal: 
\begin{equation}
\label{}
\sqrt{n}(\hat{\vartheta}_n-\vartheta) \Longrightarrow          {\cal N}(0,\II(\vartheta)^{-1}) .
\end{equation}
 Here $\II(\vartheta) $ is the Fisher information matrix 
\begin{align*}
\II(\vartheta)=\Ex_\vartheta \left[\dot \ell\left(\vartheta
  ,X_0,X_1\right)\dot \ell\left(\vartheta ,X_0,X_1\right)^\TT\right], 
\end{align*} 
where $\ell\left(\vartheta ,x,x'\right)=\ln \pi \left(\vartheta ,x,x'\right)
$. The dot means the derivation w.r.t. $\vartheta $ and $\TT$ means the 
transpose of a matrix.

As $\pi \left(\vartheta ,x,x'\right)=g\left(x'-S\left(\vartheta
,x\right)\right)$ we can write
\begin{align}
\II(\vartheta)&=\Ex_\vartheta \left[\dot \ell\left( X_j-S\left(\vartheta
  ,X_{j-1}\right)\right)\dot \ell\left( X_j-S\left(\vartheta
  ,X_{j-1}\right)\right)^\TT\right]\nonumber\\ 
&=\Ex_\vartheta \left[\frac{
    g'\left( X_j-S\left(\vartheta ,X_{j-1}\right)\right)^2\dot S\left(\vartheta
    ,X_{j-1}\right) \dot S\left(\vartheta
    ,X_{j-1}\right)^\TT}{g\left( 
    X_j-S\left(\vartheta ,X_{j-1}\right)\right)^2 }\right]\nonumber\\ 
&=\Ex\;
\left(\frac{g'\left( \varepsilon_j \right)}{g\left(\varepsilon_j \right)}
\right)^2\; \Ex_\vartheta \left[\dot S\left(\vartheta ,\xi\right) \dot
S\left(\vartheta ,\xi\right)^\TT\right]\nonumber\\ 
&=\II_g\; \Ex_\vartheta \left[\dot
  S\left(\vartheta ,\xi\right) \dot 
S\left(\vartheta ,\xi\right)^\TT\right],
\label{Fi}
\end{align}
where we  used the equality $X_j-S\left(\vartheta ,X_{j-1}\right)=\varepsilon
_j  $ and  denoted
\begin{align*}
\II_g=\int_{}^{}\frac{g'\left(x\right)^2}{g\left(x\right)}\;{\rm d}x.
\end{align*}

 Moreover the MLE is 
 asymptotically efficient. There are several definitions of the asymptotically
 efficient estimators. One of them is the following : an estimator $\vartheta
 _n^*$ is called  {\it asymptotically efficient}  if  it satisfies the relation: for all
  $\vartheta _0\in \Theta $
\begin{equation}
\label{z}
\lim_{\delta \rightarrow 0}\lim_{n \rightarrow \infty }\sup_{\left|\vartheta
  -\vartheta _0\right|<\delta }\Ex_\vartheta W\left(
\sqrt{n}\left(\vartheta _n^*-\vartheta \right)\right) = \Ex W\left(\zeta
\II(\vartheta_0)^{-1/2}\right).
\end{equation}
Here $W\left(u\right),u\in R^d$ is a loss function satisfying the usual
conditions. Note that it can be bounded, polynomial
$W\left(u\right)=\left|u\right|^p,u\in R^d$ with $p>0$ or other (see,
e.g., \cite{IH81}) and $\zeta $ is a Gaussian vector 
$\zeta \sim {\cal N}(0,{\JJ})$, ${\JJ}$ is a unit $d\times d$ matrix.
Remind that for all 
 estimators
  $\bar\vartheta _n$ the following Hajek-Le Cam's type lower bound
\begin{equation}
\label{}
\Liminf_{\delta \rightarrow 0}\Liminf_{n \rightarrow \infty
}\sup_{\left|\vartheta -\vartheta _0\right|<\delta }\Ex_\vartheta W\left(
\sqrt{n}\left(\bar\vartheta _n-\vartheta \right)\right) \geq \Ex W\left(\zeta
\II(\vartheta_0)^{-1/2}\right)
\end{equation}
holds (see, e.g. \cite{IH81}).  That is why \eqref{z} indeed defines the
asymptotically efficient estimator.

Note that these properties of the MLE were established in several works. We
mention here \cite{OI77} and  \cite{VV02} (in the one-dimensional case $d=1$).

As preliminary estimator we can use as well the BE. Recall its definition and
properties. 
Suppose that the unknown parameter
$\vartheta \in\Theta
$ is a random vector with the prior density $p\left(\vartheta  \right),
\vartheta \in \Theta $. The function $p\left(\cdot \right)$  is continuous,
bounded and positive. The BE for
the quadratic loss function  has the following representation:
\begin{align*}
\tilde\vartheta _n=\frac{\int_{\Theta }^{}\vartheta p\left(\vartheta
  \right)V\left(\vartheta ,X^n\right){\rm d}\vartheta }{\int_{\Theta }^{}
  p\left(\vartheta \right)V\left(\vartheta ,X^n\right){\rm d}\vartheta } 
\end{align*}
This estimator  under  regularity conditions is consistent, asymptotically normal   
\begin{equation}
\sqrt{n}(\hat{\vartheta}_n-\vartheta_0) \Longrightarrow N(0, \II(\vartheta_0)^{-1}) 
\end{equation}
and asymptotically efficient for the polynomial loss functions. For the proof
see \cite{MA15}. 

Recall also the properties of the estimator of the method of moments. Suppose
that the vector-function $q\left(x\right)\in R^d$ is such that the system of equations
\begin{align*}
m\left(\vartheta \right)=t,\qquad \vartheta \in \Theta 
\end{align*}
where 
\begin{align*}
m\left(\vartheta \right)=\Ex_\vartheta ^* q\left(\xi \right)
\end{align*}
has a unique solution $\vartheta =\vartheta \left(t\right)$. Introduce the
function $h\left(t\right)$ inverse to the function $m\left(\vartheta \right)$,
i.e., $\vartheta =m^{-1}\left(t\right)=h\left(t\right)$. Then the EMM is
defined as follows
\begin{align*}
\bar\vartheta _n=h\left(\frac{1}{n}\sum_{j=1}^{n}q\left(X_j\right)\right). 
\end{align*}
It is known that under regularity conditions this estimator is consistent and
asymptotically normal
\begin{align*}
\sqrt{n}\left( \bar\vartheta _n-\vartheta \right)\Longrightarrow {\cal
  N}\left(0,\CC\left(\vartheta \right) \right),
\end{align*}
where $\CC\left(\vartheta \right)$ is the  matrix defined, for example,  in
\cite{MA15}. Moreover the moments of the EMM converge too (see \cite{MA15} for
the conditions and proof). We use such estimator as preliminary one in the
numerical simulation  Example 2 below. 

In this work the construction of the multi-step MLE is based on the   score-function. Let us 
recall the definition and some properties of it.  
Introduce the log-likelihood ratio function
\begin{equation}
\label{3}
L(\vartheta,X^n)=\ln \pi(\vartheta ,X_0)+\sum_{j=1}^n \ln \pi(\vartheta,X_{j-1},X_j).
\end{equation}
The normalized   score-function  is (for simplicity of exposition we omit the
term with initial value) 
\begin{align*}
\Delta _n\left(\vartheta ,X^n\right)=\frac{1}{\sqrt{n}}\frac{\partial
  L(\vartheta,X^n)}{\partial 
  \vartheta }=\frac{1}{\sqrt{n}}\sum_{j=1}^n  \frac{g'\left(X_j-S\left(\vartheta
  ,X_{j-1}\right)\right)}{g\left(X_j-S\left(\vartheta ,X_{j-1}\right)\right)}
\dot S\left(\vartheta ,X_{j-1}\right).
\end{align*}
If we denote the true value $\vartheta =\vartheta _0$, then we have
\begin{align*}
\Delta _n\left(\vartheta_0 ,X^n\right)=\frac{1}{\sqrt{n}}\sum_{j=1}^n
\frac{g'\left(\varepsilon _j\right)}{g\left(\varepsilon _j\right)}\; \dot
S\left(\vartheta_0 ,X_{j-1}\right).
\end{align*}
 
Note that ($i<j$)
\begin{align*}
&\Ex_\vartheta  \left(\frac{g'\left(\varepsilon _i\right)}{g\left(\varepsilon _i\right)}\;\frac{g'\left(\varepsilon
  _j\right)}{g\left(\varepsilon _j\right)} \dot
S\left(\vartheta ,X_{i-1}\right)\; \dot 
S\left(\vartheta ,X_{j-1}\right)^\TT\right)\\
&\qquad =\Ex_\vartheta  \left(\left.\frac{g'\left(\varepsilon
  _i\right)}{g\left(\varepsilon _i\right)}\;\dot 
S\left(\vartheta ,X_{i-1}\right)\;\Ex_\vartheta \left(\frac{g'\left(\varepsilon
  _j\right)}{g\left(\varepsilon _j\right)}  \dot 
S\left(\vartheta ,X_{j-1}\right)^\TT\right|{\cal F}_{j-1}\right)\right)=0
\end{align*}
because
\begin{align*}
\left.\Ex_\vartheta \left(\frac{g'\left(\varepsilon
  _j\right)}{g\left(\varepsilon _j\right)}  \dot 
S\left(\vartheta ,X_{j-1}\right)^\TT\right|{\cal F}_{j-1}\right)=\Ex
\left(\frac{g'\left(\varepsilon 
  _j\right)}{g\left(\varepsilon _j\right)} \right)\left.\Ex_\vartheta \left(
\dot  
S\left(\vartheta ,X_{j-1}\right)^\TT\right|{\cal F}_{j-1}\right)
\end{align*}
and
\begin{align*}
\Ex \left(\frac{g'\left(\varepsilon
  _j\right)}{g\left(\varepsilon _j\right)} \right)=\int_{-\infty }^{\infty
}g'\left(x\right)\,{\rm d}x=0. 
\end{align*}
Therefore by the central limit theorem 
\begin{align*}
\Delta _n\left(\vartheta_0 ,X^n\right)\Longrightarrow {\cal
  N}\left(0,\II\left(\vartheta _0\right)\right),
\end{align*}
where $\II\left(\vartheta _0\right) $ is  the Fisher information matrix
defined in \eqref{Fi}.

\section{Main result}

Suppose that we have a Markov sequence $X^n=\left(X_j\right)_{j=0,\ldots,n}$ with the
transition density $\pi\left(\cdot \right) $ depending on some unknown
finite-dimensional parameter $\vartheta\in \Theta $. The set $\Theta \subset
R^d$ is open, bounded.  Our goal is to construct on-line recurrent estimator
of this parameter.  Therefore we need for each $j$ to have an estimator
$\vartheta _{j,n}^*$ with {\it good properties}, i.e., this estimator can be
easily calculated and the same time it has to be asymptotically optimal in
some sense.  We call such sequence of estimators $\vartheta _{j,n}^*,
j=1,\ldots,n$ {\it estimator-process}.  We propose a construction of such
estimator in two steps.  We slightly change the statement of the
problem. Introduce the {\it learning part}
$X^N=\left(X_0,X_1,\ldots,X_N\right)$ of observations
$X^n=\left(X_0,X_1,\ldots,X_n\right)$, where $N=\left[n^\delta\right] $ ($N $
is the integer part of $n^\delta $) and the parameter $\delta <1$ will be
chosen later.

We say that a family of random variables  $\left\{\eta _n\left(\vartheta
\right),n=1,2,\ldots\right\}$ is tight uniformly on compacts $\KK\subset\Theta $ if for any
$\varepsilon >0$ and any compact $\KK $ there exists a constant $C>0$ such
that
\begin{align*}
\sup_{\vartheta \in\KK} \Pb_\vartheta \left( \left|\eta _n\left(\vartheta
\right)\right|>C\right)\leq \varepsilon . 
\end{align*}

Throughout the paper we suppose that the following conditions
are fulfilled.

{\it Conditions} ${\cal R}$.
\begin{enumerate}
\item {\it The time series  $\left(X_j\right)_{j\geq 0} $ is strictly
  stationary and
has a unique invariant distribution with the density function $\pi 
\left(\vartheta ,x\right)$.  }

\item {\it The preliminary estimator $\bar\vartheta _n$ is such that $\sqrt{n}\left(\bar\vartheta _n-\vartheta \right) $
 is tight uniformly on compacts $\KK\subset\Theta $.} 

\item {\it The function $S\left(\vartheta ,x\right)\in {\cal C}^3_\vartheta $,
  the density $g\left(\cdot \right)>0$ and $g\left(\cdot \right)\in {\cal C}^3$.
 The derivatives ${\partial^i \ell\left(\vartheta ,x,x'\right)}/{\partial \vartheta
  ^i}, i=1,2,3 $  of the function $\ell\left(\vartheta ,x,x'\right)=\ln \pi \left(\vartheta
  ,x,x'\right)$ are uniformly on
  $\vartheta $  majorated by quadratically integrable functions, i.e., }
\begin{align*}
\sup_{\vartheta \in \Theta }\left\|\frac{\partial^i \ell\left(\vartheta
  ,x,x'\right)}{\partial \vartheta 
  ^i}\right\|\leq R_i\left(x,x'\right),\qquad i=1,2,3, 
\end{align*}
{\it where  $\Ex _\vartheta \left| R_i\left(X_{j-1},X_{j}\right) \right|^2<C$
and the constant $C>0$ does not depend on $\vartheta $.}
\item {\it  We have }
\begin{itemize}
\item {\it the law of large numbers
\begin{align}
\label{lln}
&\frac{1}{n}\sum_{j=1}^{n}\dot \ell\left(\vartheta ,X_{j-1},X_j\right) \dot
\ell\left(\vartheta ,X_{j-1},X_j\right)^\TT  \longrightarrow
\II\left(\vartheta \right),
\end{align}
\item the central limit theorem 
\begin{align}
\label{clt1}
&\frac{1}{\sqrt{n}}\sum_{j=1}^{n}\dot \ell\left(\vartheta
,X_{j-1},X_j\right) \Longrightarrow {\cal N}\left(0,\II\left(\vartheta
\right) \right),
\end{align}
\item the family of random variables 
\begin{align}
\label{clt2}
&\frac{1}{\sqrt{n}}\sum_{j=1}^{n}\left[\ddot \ell\left(\vartheta
,X_{j-1},X_j\right)+
\II\left(\vartheta \right)\right] 
\end{align}
is tight  uniformly on compacts $ \KK\subset \Theta $. }
\end{itemize}
\item {\it  The information  matrix $\II\left(\vartheta \right)$  is
  Lipschtitz
\begin{align}
\label{Lip}
\left|\II\left(\vartheta_1 \right)-\II\left(\vartheta_2 \right)\right|\leq
L\left|\vartheta_1 -\vartheta_2 \right|
\end{align}
 and is uniformly
  in $\vartheta \in \Theta $ non-degenerate and bounded 
\begin{equation}
\label{FM}
0<\inf_{\vartheta \in \Theta }\inf_{\left|\lambda \right|=1}\lambda ^\TT
\II(\vartheta)\lambda ,\qquad \sup_{\vartheta \in \Theta }\sup_{\left|\lambda
  \right|=1}\lambda ^\TT \II(\vartheta)\lambda <\infty . 
\end{equation}
Here $\lambda \in R^d$. }
\end{enumerate}

Note that as preliminary estimator $\bar\vartheta _N$ we can take the MLE, the
BE or the EMM. All of them have the required properties (under additional
regularity conditions, which we do not mention here). The details can be found
in  \cite{OI77},
\cite{VV02}, \cite{MA15} or any other work describing their properties.  
  The conditions for \eqref{lln}-\eqref{clt2} can be found, for example, in
\cite{FY03},  \cite{TK00},\cite{Ver00}. The condition \eqref{Lip} can be verified if we
have  the corresponding smoothness of the density of invariant distribution
$\pi \left(\vartheta ,x\right)$ (see, e.g. \cite{BSV15}).

 We construct the one-step MLE-process $ \vartheta
_{k,n}^\star,k=N+1,\ldots,n$ as follows. Introduce the variable $s\in
\left[\tau _\delta ,1\right]$, where $\tau _\delta =n^{-1+\delta }\rightarrow 0$ and put
$k=\left[sn\right]$, where $\left[a\right]$ means the integer part of $a$. Let
us write $\vartheta _{k,n}^\star=\vartheta _{s,n}^\star$ and consider the
estimator-process  $\vartheta ^\star_n=\left(\vartheta _{s,n}^\star, s\in
\left[\tau _\delta ,1\right]\right)$. Our 
goal is to construct an estimator process $\vartheta ^\star_n $ asymptotically
optimal for all $s\in \left[\tau _\delta ,1\right]$. Recall that the MLE
$\hat\vartheta _{s,n}$ constructed by the first $k=\left[sn\right] $
observations is asymptotically efficient and for example, 
$$
\sqrt{sn}\left(\hat\vartheta _{s,n}-\theta \right)\Longrightarrow {\cal
  N}\left(0,\II\left(\vartheta \right)^{-1}\right), \quad s\in \left[\delta ,1\right].
$$
Note that to solve the  equation
\begin{align*}
\sup_{\vartheta \in \Theta }V\left(\vartheta
,X^{\left[sn\right]}\right)=V\left(\hat\vartheta_{s,n}
,X^{\left[sn\right]}\right) 
\end{align*}
for all $s\in \left[\tau _\delta ,1\right]$ is computationally rather difficult
problem, except some particular examples. Therefore it is better to seek another
estimators, which have the same limit covariance matrix as the MLE (which is
asymptotically efficient) for 
all $s\in \left(\tau _\delta ,1]\right)$ and which can be calculated in more
  simple way. 

We consider two different situations depending on the length of the learning
interval $\left[0,N\right]$. If $N=\left[n^\delta\right] $ (here
$\left[a\right] $ is integer part of $a$) with $\frac12<\delta <1 $ then
we construct the one-step MLE-process and if we take the preliminary interval
shorter, i.e., $N=\left[n^\delta\right] $ with $\frac14<\delta \leq \frac12 $, then we
introduce an intermediate estimator and only after that we can construct the
two-step MLE-process. Therefore we consider below these two situations separately.

\subsection{Case $N=\left[n^{\delta} \right],\frac12<\delta < 1$}

We proceed as follows. 
Let us fix $s\in \left[\tau _\delta ,1\right]$ and slightly modify the vector score-function  
\begin{align*}
\Delta_k(\vartheta,X^k_{N})=\frac{1}{\sqrt{k}} \sum_{j=N+1}^k \dot \ell
(\vartheta, X_{j-1}, X_j) ,
\end{align*}
where  $k=\left[sn\right]\rightarrow \infty $.  Introduce the one-step MLE-process
\begin{align*}
\vartheta_{s,n}^\star=\bar \vartheta _N + \frac{1}{\sqrt k}\II\left(\bar \vartheta
_N\right)^{-1}{\Delta_k(\bar \vartheta _N,X^k_{N})}       { } ,\qquad \tau _\delta \leq s\leq 1
\end{align*}

 Here and below for simplicity of
notation this writing
means that $N$ is the integer part of $n^{\delta} $.

\begin{theorem}
\label{T11} Suppose that the conditions ${\cal R}$  are fulfilled, then for
all $s\in (0,1]$
\begin{align}
\label{an}
\sqrt{k}(\vartheta_{s,n}^\star-\vartheta)\Longrightarrow {\cal
  N}\left(0,\II\left(\vartheta \right)^{-1}\right)
\end{align}
and this estimator-process is asymptotically efficient for the bounded loss
functions in \eqref{z}. 
\end{theorem}
{\bf Proof.} Note that for any $s>0$ ($s\leq 1$) we have $s>\tau _\delta $ for
$n>s^{\frac{1}{1-\delta }} $.    We can write
\begin{align*}
\sqrt{k}(\vartheta_{s,n}^{\star}-\vartheta)&=\sqrt{k}(\overline{\vartheta}_N-\vartheta)+\II\left(\bar
\vartheta _N\right)^{-1}{\Delta_k(\bar \vartheta _N,X^k_{N})}\\
&=\sqrt{k}(\overline{\vartheta}_N-\vartheta)+\II(\bar{\vartheta}_N
{}{)}^{-1}\Delta_k(\vartheta,X^k_{N})\\
&\qquad +\II(\bar{\vartheta}_N
{}{)}^{-1}\left[\Delta_k(\bar\vartheta,X^k_{N})-\Delta_k(\vartheta,X^k_{N}) \right].
 \end{align*}
We have
\begin{align*}
\Delta_k(\bar\vartheta,X^k_{N})-\Delta_k(\vartheta,X^k_{N})=\int_{0}^{1}\langle
(\bar{\vartheta}_N-\vartheta), \dot
\Delta_k(\vartheta+v\left(\bar{\vartheta}_N-\vartheta \right),X^k_{N})\rangle\,{\rm d}v .
\end{align*}
Hence (below $\vartheta _v=\vartheta+v\left(\bar{\vartheta}_N-\vartheta \right)$)
\begin{align*}
&\sqrt{k}(\overline{\vartheta}_N-\vartheta)+\II(\bar{\vartheta}_N
)^{-1}\left[\Delta_k(\bar\vartheta,X^k_{N})-\Delta_k(\vartheta,X^k_{N})
  \right]\\
&\qquad =\sqrt{k}(\overline{\vartheta}_N-\vartheta)\II(\bar{\vartheta}_N
)^{-1}\left[ \II(\bar{\vartheta}_N
)+ \frac{1}{\sqrt{k}}\int_{0}^{1}
 \dot
\Delta_k(\vartheta _v,X^k_{N})\,{\rm d}v    \right].
\end{align*}
Further
\begin{align*}
&\II(\bar{\vartheta}_N )+ \frac{1}{\sqrt{k}}\int_{0}^{1} \dot
\Delta_k(\vartheta _v,X^k_{N})\,{\rm d}v=\II({\vartheta})+ \frac{1}{\sqrt{k}} \dot
\Delta_k(\vartheta ,X^k_{0})-\frac{1}{\sqrt{k}} \dot
\Delta_k(\vartheta ,X^{N-1}_{0})\\
&\qquad +\II(\bar{\vartheta}_N )-\II({\vartheta})+\frac{1}{\sqrt{k}}\int_{0}^{1} \left[\dot
\Delta_k(\vartheta _v,X^k_{N})- \dot
\Delta_k(\vartheta ,X^k_{N}) \right]\,{\rm d}v\\
&\quad =\frac{1}{k}\sum_{j=1}^{k}\left[\ddot \ell\left(\vartheta
  ,X_{j-1},X_j\right)+
  \II({\vartheta})\right]+O\left(\frac{N}{k}\right)+O\left( n^{-\frac{\delta }{2}}\right),
\end{align*}
because
\begin{align*}
&\frac{1}{\sqrt{k}} \dot \Delta_k(\vartheta ,X^{N-1}_{0})
  =\frac{1}{k}\sum_{j=1}^{N-1}\ddot\ell \left(\vartheta
  ,X_{j-1},X_j\right)=O\left(\frac{N}{k}\right)=O\left( n^{-1+\delta }\right),\\
&\left| \II(\bar{\vartheta}_N
)-\II({\vartheta})\right| \leq L\left| \bar{\vartheta}_N-\vartheta\right|
=O{\left(n^{-\frac{\delta }{2}}\right)}
\end{align*}
and
\begin{align*}
\frac{1}{\sqrt{k}}\int_{0}^{1} \left[\dot
\Delta_k(\vartheta _v,X^k_{N})- \dot
\Delta_k(\vartheta ,X^k_{N}) \right]\,{\rm d}v=O\left(
\bar{\vartheta}_N-\vartheta \right)=O{\left(n^{-\frac{\delta }{2}}\right)}. 
\end{align*}
Here and in the sequel $O\left(n^{-c}\right)$ means  that
$n^{c}O\left(n^{-c}\right) $ is bounded in probability uniformly on compacts $\KK$, i.e., for any
$\varepsilon >0$ there exists $C_1>0$ such that
\begin{align*}
\sup_{\vartheta \in \KK}\Pb_\vartheta \left( n^{c}\left|O\left(n^{-c}\right)\right|>C_1  \right)\leq
\varepsilon . 
\end{align*}
For example,
\begin{align*}
&\sup_{\vartheta \in \KK}\Pb_\vartheta \left(  \frac{1}{k}\sum_{j=1}^{N-1}\left|\ddot\ell \left(\vartheta
  ,X_{j-1},X_j\right) \right|   >C_1  \right)\leq \frac{1}{kC_1}
\sum_{j=1}^{N-1}\sup_{\vartheta \in \KK}\Ex_\vartheta  \left|\ddot\ell \left(\vartheta
  ,X_{j-1},X_j\right) \right|\\
&\qquad \quad \leq \frac{1}{kC_1}
\sum_{j=1}^{N-1}\sup_{\vartheta \in \KK}\Ex_\vartheta \left|R_2\left( X_{j-1},X_j\right)  \right|\leq \frac{CN}{C_1k}.
\end{align*}

Recall that $\Ex_\vartheta
\ddot\ell(\vartheta,X_{j-1},X_j)=-\II({\vartheta}). $ Hence by the central
limit theorem \eqref{clt2} we have 
\begin{align*}
\frac{1}{{\sqrt{k}}}\sum_{j=1}^{k}\left[\ddot\ell(\vartheta,X_{j-1},X_j)
  +\II({\vartheta} 
{}{)}\right]\Longrightarrow {\cal N}\left(0,\DD\left(\vartheta \right)\right) 
\end{align*} 
with some  $\DD\left(\vartheta \right)$. 

Therefore 
\begin{align*}
&\sqrt{k}(\vartheta_{s,n}^{\star}-\vartheta)=\II(\bar{\vartheta}_N
  {}{)}^{-1}\Delta_k(\vartheta,X^k_{N})\\ &\qquad +n^{\frac{\delta
    }{2}}(\overline{\vartheta}_N-\vartheta) \left[n^{\frac{1-\delta }{2}}
    O\left(n^{-\frac{1}{2}}\right)+n^{\frac{1-\delta }{2}} O\left(n^{-1+\delta
    }\right)+n^{\frac{1-\delta }{2}}O\left(n^{-\frac{\delta }{2}} \right)
    \right]\\ &\quad =\II({\vartheta}
  {}{)}^{-1}\Delta_k(\vartheta,X^k_{0})+o\left(1\right)\Longrightarrow {\cal
    N}\left(0,\II({\vartheta} {}{)}^{-1}\right),
\end{align*}
where we used once more the central limit theorem \eqref{clt1}.

Therefore the one-step MLE-process $\vartheta _n^\star =\left(\vartheta
_{s,n}^\star,\tau _\delta <s \leq 1\right)$ for all $s\in  (\tau _\delta,  1]$
  is uniformly in $\vartheta \in \KK$ asymptotically normal \eqref{an}. Hence
  for the bounded loss functions $W\left(\cdot \right)$
  we obtain the convergence
\begin{align*}
\lim_{n \rightarrow \infty }\sup_{\left|\vartheta
  -\vartheta _0\right|<\delta }\Ex_\vartheta W\left(
\sqrt{n}\left(\vartheta _n^*-\vartheta \right)\right) = \sup_{\left|\vartheta
  -\vartheta _0\right|<\delta }\Ex W\left(\zeta
\II(\vartheta)^{-1/2}\right).
\end{align*}
Now \eqref{z} follows from the continuity of the Fisher information.

\subsection{Case $N=n^{\delta} ,\frac14<\delta \leq \frac12$}

The choice of the  learning period of observations $N=\left[n^\delta \right]$
with $\delta \in (1/2,1)$
allows us to construct an estimator process for the values $s\in (\tau _\delta
,1]  $ only. It can be interesting to see if it is possible to take more short
  learning interval and therefore to have the estimator-process for the larger
  time interval. Our goal is to show that the learning period can be  $N=\left[n^\delta \right]$
with $\delta \in (1/4,1/2]$. Below we follow the construction which
was already realized in \cite{Kut15} in the case of ergodic diffusion
process. 

Suppose that $N=\left[n^\delta \right]$
with $\delta \in (1/4,1/2]$. The asymptotically efficient estimator we construct in three steps. By the
first $N$ observations as before  we obtain the preliminary estimator $\bar\vartheta _{N}$ 
which is asymptotically normal  with the rate $\sqrt{N}$, i.e.,
$$
n^{\frac{\delta }{2}}\left(\bar\vartheta _{N}-\vartheta \right) \Longrightarrow {\cal
  N}\left(0, \BB\left(\vartheta \right)\right). 
$$
This can be the same estimator as in the preceding case.  It can be, for
example, the EMM, BE or MLE.

The two-step MLE-process $\vartheta
_n^{\star \star}=\left(\vartheta _{s,n}^{\star \star},k=N+1,\ldots,n\right)$
we construct as follows.  
Fix some $s\in (\tau _\delta ,1]$, $\tau _\delta
=n^{-1+\delta }$ and introduce the second preliminary estimator-process
(as before  $k=\left[sn\right]$)
\begin{align}
\label{kk}
\bar\vartheta_{k,2}=\bar \vartheta _{N} + \frac{1}{\sqrt k}\II\left(\bar \vartheta
_{N}\right)^{-1}{\Delta_k(\bar \vartheta _{N},X^k)},\qquad k=N+1,\ldots,n
,\end{align}
where
\begin{align*}
\Delta_k(\vartheta ,X^k)=\frac{1}{\sqrt{k}}\sum_{j=1}^{k}\dot \ell\left(\vartheta ,X_{j-1},X_j\right).
\end{align*}
Then  we  show that the random sequence ${n^{1/4+\varepsilon  }}\left(\bar\vartheta _{k,2}-\vartheta
\right)  $ with some $\varepsilon >0$ is bounded in probability (tight). 

Finally, using this estimator-process and the one-step procedure of  Theorem
\ref{T11} we obtain the asymptotically efficient estimator
\begin{align}
\label{ae}
\vartheta_{s,n}^{\star \star}=\bar \vartheta _{k,2} + \frac{1}{\sqrt k}\II\left(\bar \vartheta
_{k,2}\right)^{-1}{\Delta_k(\bar \vartheta _{k,2},X^k)}  .
\end{align}
In the next theorem we realize this program. 

\begin{theorem}
\label{T22} Suppose that the conditions of regularity are fulfilled, then
the estimator $\vartheta _{s,n}^{\star \star}$ defined  by \eqref{kk} and
\eqref{ae} for all $s\in (0,1]$ is asymptotically normal
\begin{align*}
\sqrt{k}(\vartheta_{s,n}^{\star \star}-\vartheta)\Longrightarrow {\cal
  N}\left(0,\II\left(\vartheta \right)^{-1}\right)
\end{align*}
and asymptotically efficient for the bounded loss functions.
\end{theorem}
{\bf Proof.} The only thing to proof is the tightness of the sequence of
random vectors ${n^{1/4+\varepsilon }}\left(\bar\vartheta _{k,2}-\vartheta
\right) $  because if it is tight, then the proof of Theorem \ref{T22} follows from the Theorem
\ref{T11}. Let us fix some $\varepsilon \in \left(0,\frac14\right)$. 

For the estimator-process $\bar\vartheta
_{k,2} $ defined by \eqref{kk} we can write
\begin{align*}
&{n^{\frac{1}{4}+\varepsilon }}\left(\bar\vartheta _{k,2}-\vartheta\right)
=n^{\frac{1}{4}+\varepsilon }\left(\bar\vartheta _{N}-\vartheta\right)+
\frac{n^{\frac{1}{4}+\varepsilon }}{\sqrt k}\II\left(\bar \vartheta
_{N}\right)^{-1}{\Delta_k(\bar \vartheta _{N},X^k)}\\
&\qquad =n^{\frac{1}{4}+\varepsilon  }\left(\bar\vartheta
_{N}-\vartheta\right)+ 
\frac{n^{\frac{1}{4}+\varepsilon }}{\sqrt k}\II\left(\bar \vartheta
_{N}\right)^{-1}{\Delta_k(\vartheta,X^k)}\\
&\qquad\quad  + 
\frac{n^{\frac{1}{4}+\varepsilon }}{\sqrt k}\II\left(\bar \vartheta
_{N}\right)^{-1}{\left(\bar \vartheta
_{N}-\vartheta  \right)  \dot\Delta_k(\tilde \vartheta _{N},X^k)}.
\end{align*}

Note that $\Delta_k(\vartheta,X^k) $ is asymptotically normal and therefore 
\begin{align*}
\frac{n^{\frac{1}{4}+\varepsilon }}{\sqrt k}\II\left(\bar \vartheta
_{N}\right)^{-1}{\Delta_k(\vartheta,X^k)}\longrightarrow 0,
\end{align*}
because ${n^{\frac{1}{4}+\varepsilon }}{ k^{-\frac{1}{2}}}\rightarrow 0
$. Further
\begin{align*}
&n^{\frac{1}{4}+\varepsilon }\left(\bar\vartheta
_{N}-\vartheta\right) + 
\frac{n^{\frac{1}{4}+\varepsilon }}{\sqrt k}\II\left(\bar \vartheta
_{N}\right)^{-1}{\left(\bar \vartheta
_{N}-\vartheta  \right)  \dot\Delta_k(\tilde \vartheta _{N},X^k)}\\
&\qquad =n^{\frac{1}{8}+\frac{\delta }{2} }\left(\bar\vartheta_{N}-\vartheta\right)R_n,
\end{align*}
where
\begin{align*}
R_n=n^{\frac{1}{8}+\varepsilon -\frac{\delta }{2} }\left[\JJ +\II\left(\bar \vartheta
_{N}\right)^{-1} \frac{1}{k}\sum_{j=1}^{k} \ddot \ell\left(\tilde\vartheta_N
  ,X_{j-1},X_j\right)\right]. 
\end{align*}
We have by the law of large numbers
\begin{align*}
 \frac{1}{k}\sum_{j=1}^{k} \ddot \ell\left(\vartheta
 ,X_{j-1},X_j\right)\longrightarrow -\II\left( \vartheta \right).
\end{align*}
From the regularity conditions it follows that
\begin{align*}
&\left|\II\left(\bar \vartheta
_{N}\right)^{-1}-\II\left( \vartheta \right)^{-1}\right|\leq C\left|\bar \vartheta
_{N}-\vartheta  \right|,\\
&k^{-1/2} \left|\dot\Delta_k(\tilde \vartheta _{N},X^k)-\dot\Delta_k( \vartheta
  _{k},X^k)\right|\leq \frac{1}{k}
  \sum_{j=1}^{k}\left|R_3\left(X_{j-1},X_j\right)\right|\left|\bar \vartheta
  _{N}-\vartheta  \right|. 
\end{align*}
Therefore we verified the tightness of the sequence
${n^{\frac{1}{4}+\varepsilon }}\left(\bar\vartheta _{k,2}-\vartheta\right)
$. Now the proof of  the Theorem \ref{T22} follows from the proof of the Theorem \ref{T11}.

\section{Examples} 

We consider below two examples. The first one is new and the second example
was already discussed in the previous work in the context of the study of the
Bayesian estimators and the estimators of the method of moments
\cite{MA15}. In the first 
example we construct the preliminary MLE  and  the one-step
MLE-process. In the second example we construct the preliminary EMM,  the
second preliminary estimator-process and then the two-step MLE-process.  

\subsection{Example 1.}
Let us consider the problem of the construction of the one-step MLE-process in
the case of observations $X^n=\left(X_0,X_1,\ldots,X_n\right)$  of the time series
\begin{align}
X_{j}&=\frac{(X_{j-1})^2 }{1+\vartheta \mid{X_{j-1}}\mid}+\varepsilon _{j},
\,\qquad \vartheta \in (2, 5) ,
\label{equat11}
\end{align}
where $ \left(\varepsilon _{j}\right)_{j\geq 1}\sim {\cal N}(0,1)$. 

Note that  this time series has invariant distribution. 
  The
density of it  we estimate with the help  of gaussian kernel-type
estimator $ K\left(\cdot \right)$: 
\begin{align*}
\hat{\pi }_n(x)=\frac{1}{n h_n}\sum_{j=1}^n K \left(\frac{X_j-x}{h_n}\right),\qquad 
K(x)=\frac{1}{\sqrt{2 \pi}}\;e^{-\frac{x^2}{2}} , 
\end{align*}
where the width $h_n=n^{-1/5} $.

\begin{figure}[h]
\centering \includegraphics [scale=0.50]  {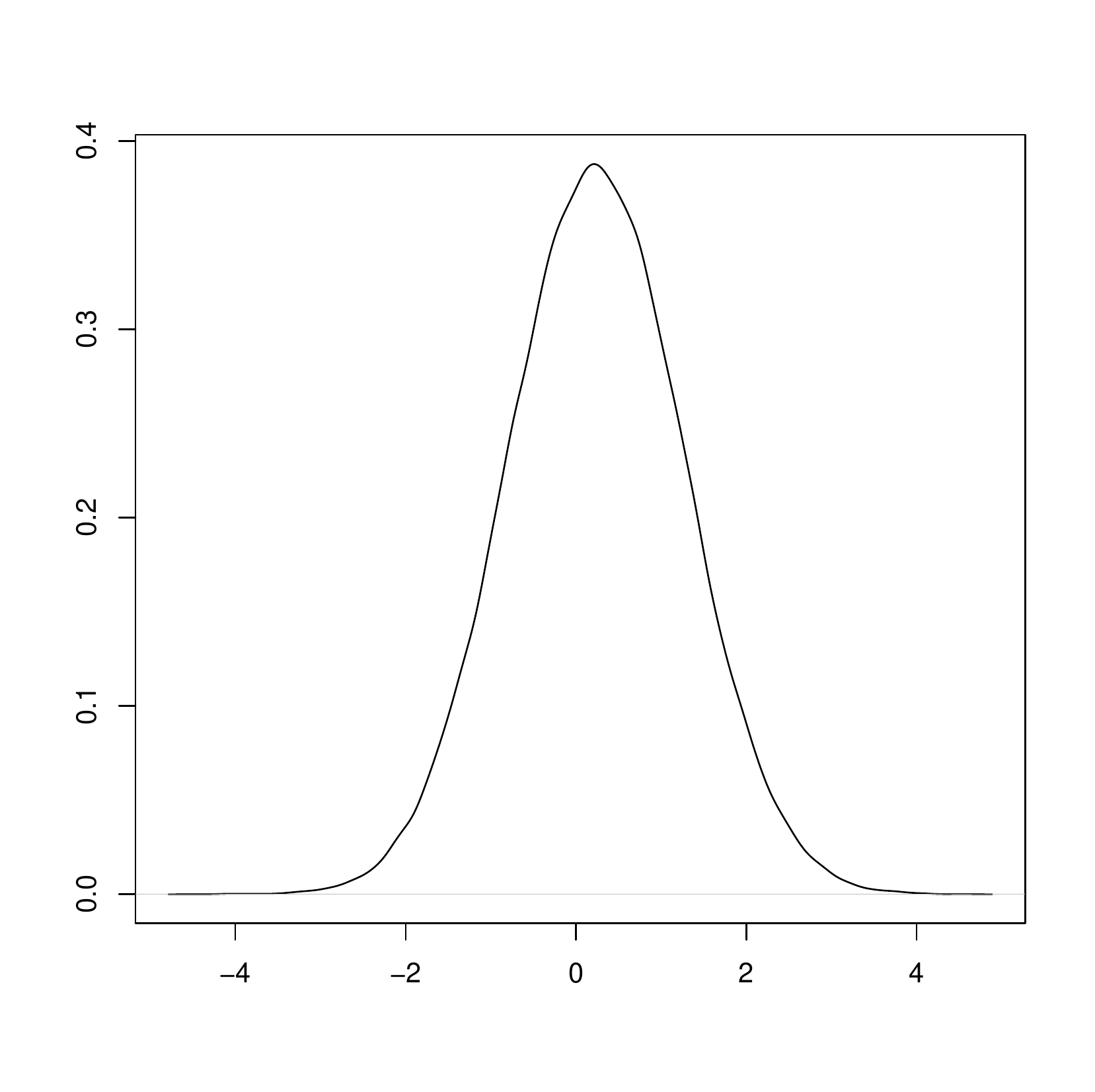}
\caption{Estimator of invariant density $\pi^\star(\vartheta, x)$ for
$\vartheta=2.5$ and $n=10^5$} 
\label{figure_30}
\end{figure}

On the Figure \ref{figure_30} we present the estimator of the invariant
density in the case  $n=10^5$ and $\vartheta =2,5$. 

First we define the MLE constructed on the learning sequence
$X^N=\left(X_0,X_1,\ldots,X_N\right)$.  For the conditional density function $\pi (\vartheta, X_{j-1}, X_{j})$ of
the Markov sequence \eqref{equat11}, we have the representation
\begin{align}
\pi (\vartheta, X_{j-1}, X_{j})=\frac{1}{\sqrt{2\pi}}e^{-\frac{1}{2} \left[X_j-\frac{(X_{j-1})^2 }{1+\vartheta \mid{X_{j-1}}\mid}\right]^2}.
\label{eq22}
\end{align}
Hence  the log-likelihood ratio function is 
\begin{align*}
L_N(\vartheta,X^N)&=\ln \pi _0\left(X_0\right)+\sum_{j=1}^N \left(-\frac{1}{2}
\ln {2 \pi}-\frac{1}{2} \left[X_j-\frac{(X_{j-1})^2 }{1+\vartheta
    \mid{X_{j-1}}\mid}\right]^2 \right)\\ 
&=\ln \pi _0\left(X_0\right)+\sum_{j=1}^N \ell (\vartheta,X_{j-1},X_j),\qquad \vartheta\in (2, 5).
\end{align*}

 To find the MLE $\hat\vartheta _{N}$ we have to
solve the maximum likelihood equation 
\begin{align*}
\frac{\partial L}{\partial \vartheta}= \sum_{j=1}^{N}\dot\ell\left(\vartheta,
X_{j-1}, X_j\right)= 0, \qquad \vartheta \in \left(2, 5\right),
\end{align*}
 which has the following form
\begin{align*}
\sum_{j=1}^N  \frac{\mid{X_{j-1}}\mid^3 }{(1+\vartheta
  \mid{X_{j-1}}\mid)^2}\left(-X_j+\frac{(X_{j-1})^2 }{1+\vartheta
  \mid{X_{j-1}}\mid}\right)=0,\qquad \vartheta \in \left(2, 5\right).
\end{align*}


Now we construct the one-step MLE-process  $\vartheta
_n^\star=\left(\vartheta _{k,n}^\star,N+1\leq k\leq n\right)$ based 
on this preliminary estimator $\hat \vartheta _N$ as follows. The normalized 
score-function is 
\begin{align*}
\Delta_k(\vartheta ,X^k)=\frac{1}{\sqrt{k}}\sum_{j=1}^{k}\frac{\mid{X_{j-1}}\mid^3 }{(1+\vartheta
  \mid{X_{j-1}}\mid)^2}\left(-X_j+\frac{(X_{j-1})^2 }{1+\vartheta
  \mid{X_{j-1}}\mid}\right),
\end{align*}
where    $ N+1\leq k\leq n$. Finally the one-step MLE-process has the following representation
\begin{align*}
\vartheta_{k,n}^\star=\hat \vartheta _{N} + \frac{1}{\II_k(\hat\vartheta _N)
  {k}}\sum_{j=1}^{k}\frac{\mid{X_{j-1}}\mid^3 }{(1+\hat\vartheta _{N}
  \mid{X_{j-1}}\mid)^2}\left(-X_j+\frac{(X_{j-1})^2 }{1+\hat \vartheta _{N}
  \mid{X_{j-1}}\mid}\right)   ,
\end{align*}
where    $N+1\leq k\leq n $ and $\II_k(\hat\vartheta _N)=0.001$ is the Fisher information calculated as follows
\begin{align*}
\II_k(\hat\vartheta _N)=-\frac{1}{k}\sum_{j=1}^{k}\ddot\ell
\left(\hat\vartheta _N,X_{j-1},X_j\right) .
\end{align*}
More detailed analysis shows that with such definition of the empirical Fisher
information the main result of this work Theorem 1 is valid. Therefore the
estimator-process $ \vartheta _n^\star$ is asymptotically normal with the same
limit variance as that of the MLE.  

\begin{figure}[h]
\centering\includegraphics [scale=0.68]  {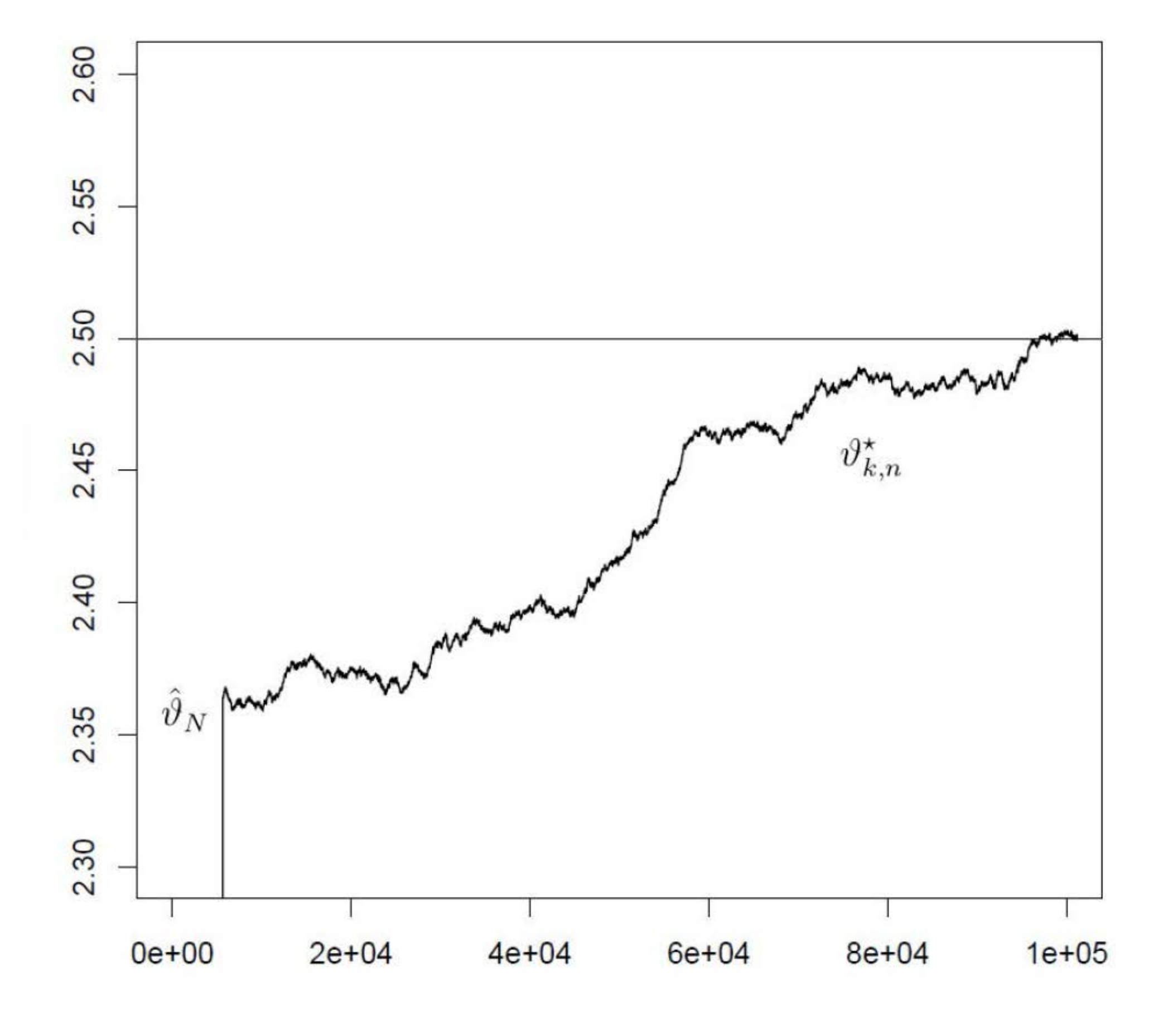}
\caption{One-step MLE-process for $n=10^5$ observations and $\vartheta=2.5$}
\label{figure_32}
\end{figure}

The  realization of the simulated one-step MLE-process for $n=10^5$ is shown on the Figure
\ref{figure_32}. We can see that the initial estimator $\hat\vartheta _{N} $
is far from the true value and that the   trajectory of one-step MLE-process
approaches to the true value. 

\subsection{Example 2.}

Let us consider another example, where preliminary estimator is 
EMM. Our goal is to illustrate the convergence of the one- and 
two-step MLE-processes, when the initial estimator is EMM (``bad'' rate and
``bad'' limit variance). 

 Introduce the time series 
\begin{align}
\label{exMain}
X_{j}=X_{j-1}+3\frac{\vartheta -X_{j-1}}{1+\left(X_{j-1}-\vartheta \right)^2}+\varepsilon
_{j}, \quad j=1,\ldots,n,
\end{align}
where $\left(\varepsilon _j\right)_{j\geq 1}$ are i.i.d. standard Gaussian
random variables and $X_0$ is given. The unknown parameter $\vartheta \in \Theta
=\left(-1,1\right)$.  This
example was already discussed in the  work \cite{MA15} to illustrate the
properties of  the BE and EMM. 

This process has ergodic properties and its invariant density can be estimated
as in the Example 1 with the help of the kernel-type estimator. The result of
such estimation can be found in \cite{MA15}.

We construct two estimator-processes: one-step and two-step.  Our goal is to
construct the estimator-processes $\vartheta _n^\star$ and $\vartheta
_n^{\star \star}$ , which are asymptotically equivalent to the MLE and
therefore are asymptotically efficient. The same time their calculation is
much more simple than that of the MLE.

We start with the one-step MLE-process. As described before we construct this
estimator in two steps.  
First we need to calculate a consistent preliminary estimator
$\bar\vartheta _N$ by the initial observations $X_1,\ldots,X_N$, where
$N=n^\delta $ with $\delta \in (\frac{1}{2},1)$.
   Note that the unknown
 parameter for this model of observations is the shift parameter and that the
 invariant density function is symmetric with respect to $\vartheta $. Hence
 we can take the EMM
\begin{align*}
\bar\vartheta_N=\frac{1}{N}\sum_{j=1}^{N} X_j\longrightarrow \vartheta , \qquad N=n^{3/4}.
\end{align*}
Of course, the limit variance of the EMM $\bar\vartheta_N $ is greater than
that of the BE, but this estimator is much more easier to calculate.

\begin{figure}[h]
\centering\includegraphics [scale=0.73]  {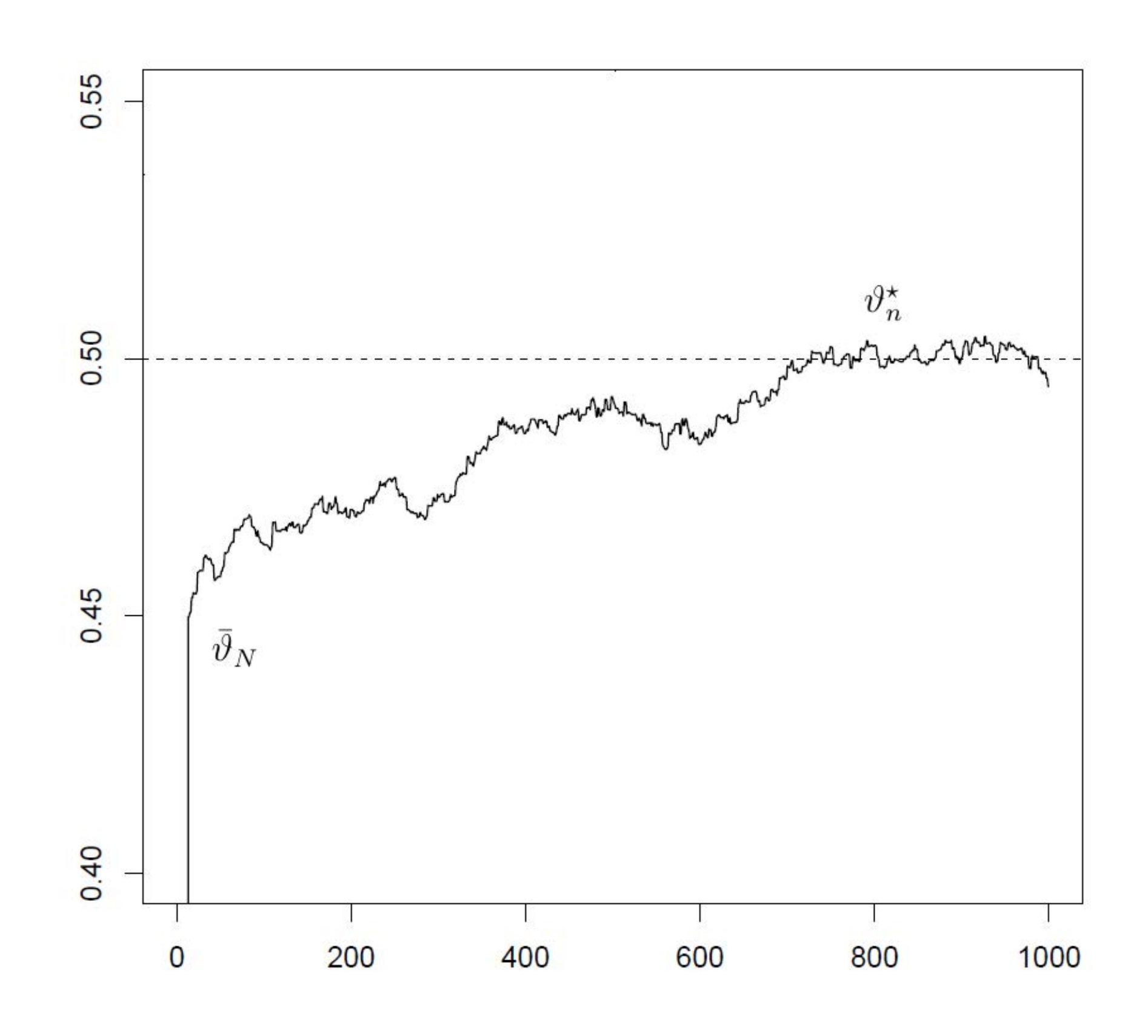}
\caption{One-step MLE-process for $n=1 000$ and $\vartheta=0.5$}
\label{figure_20}
\end{figure}

 The  score-function process is 
\begin{align*}
\Delta_k(\vartheta
,X^k)=\frac{1}{\sqrt{k}}\sum_{j=1}^{k}\dot\ell\left(\vartheta, X_{j-1},
X_j\right) ,\qquad N+1\leq k\leq n.
\end{align*}
where
\begin{align*}
\dot\ell\left(\vartheta, x, x'\right)=
&3\, \left(x'-x-3\frac{\vartheta -x}{1+\left(\vartheta-x \right)^2}\right)\, \frac{1-(\vartheta -x)^2}{\left(1+(\vartheta-x )^2\right)^2}.
\end{align*}
Therefore we can calculate the one-step MLE-process as follows
\begin{align*}
\vartheta_{k,n}^\star&=\bar \vartheta _{N} \\
&\quad + \frac{3}{\II_k
  {k}}\sum_{j=1}^{k}\left(X_j-X_{j-1}-3\frac{\bar \vartheta _{N}
  -X_{j-1}}{1+\left(\bar \vartheta _{N}-X_{j-1} \right)^2}\right)\, \frac{1-(\bar \vartheta _{N}
  -X_{j-1})^2}{\left(1+(\bar \vartheta _{N}-X_{j-1}  )^2\right)^2} .
\end{align*}
Here $\II_k $ is the empirical Fisher information. Its calculation in this
example can be found in \cite{MA15}. Note that $\II\left(\vartheta
\right)=\II$ as usual with the shift parameter. 
Remind that by the Theorem \ref{T11}  this estimator  is asymptotically normal.  

\begin{figure}[h]
\centering\includegraphics [scale=0.73]  {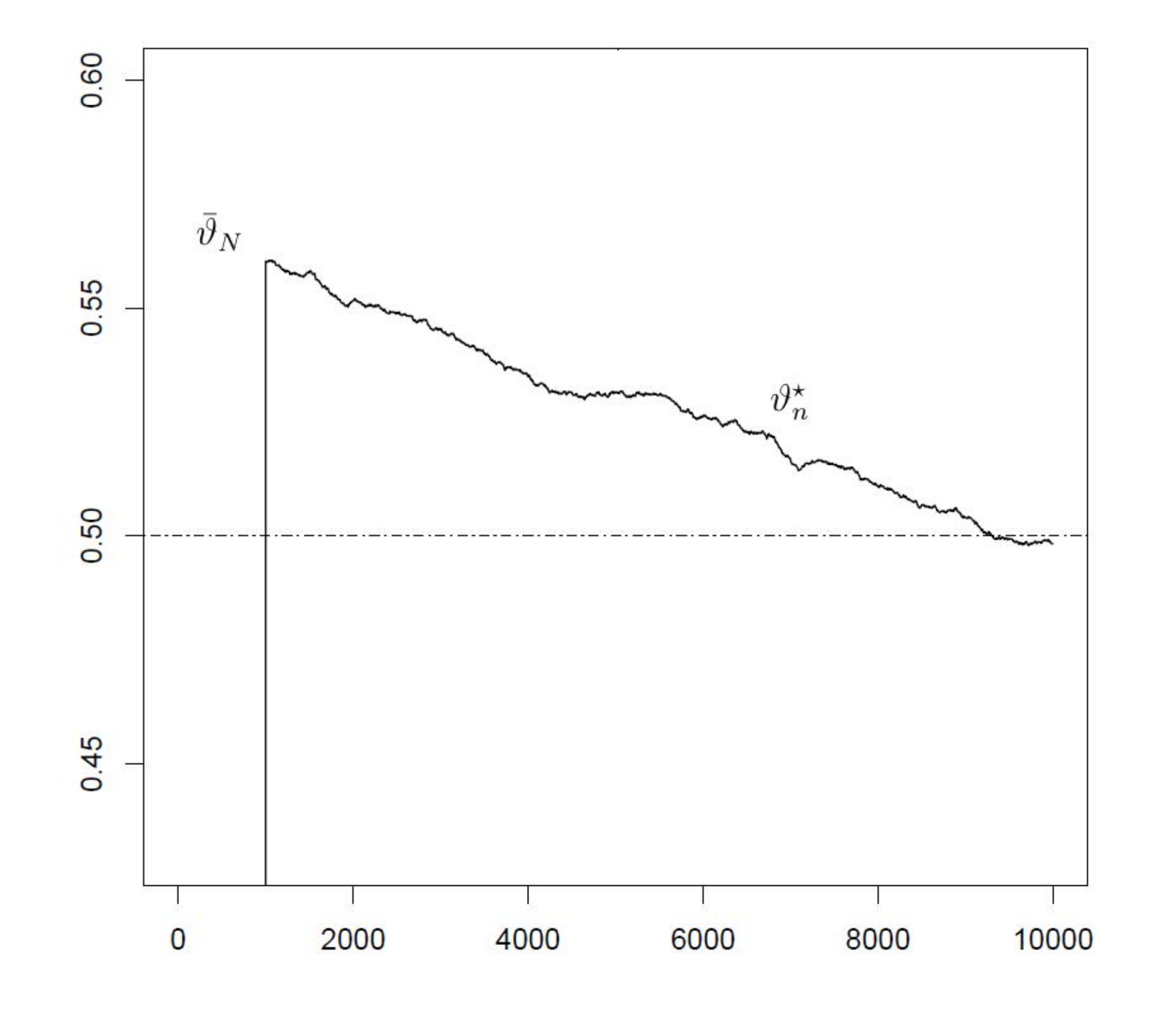}
\caption{One-step MLE-process for $n=10^4$ and $\vartheta=0.5$}
\label{figure_21}
\end{figure}

The simulated one-step MLE-processes are shown on the Figure \ref{figure_20}
and \ref{figure_21} for $n=10^3$ and $n=10^4$ respectively.

On the Figure \ref{figure_20} the preliminary EMM $\bar\vartheta _N=0.45$ that
is  close to the true value of parameter $\vartheta =0.5$. We obtain this
estimator based on the learning interval of $N=178$ observations. And we can
observe the sequence of estimator $\vartheta _n^\star=\left(\vartheta
_{k,n}^\star, k=N+1;\ldots,n\right)$ that is asymptotically efficient.

On the Figure \ref{figure_21} the preliminary EMM $\bar\vartheta _N=0.56$ that
is  close to the true value $\vartheta =0.5$. We obtain this estimator
by the first  $N=10^3$ observations.  We can see that 
the  estimator-process  $\vartheta _n^\star=\left(\vartheta _{k,n}^\star, 
k=N+1;\ldots,n\right)$  tends to the true value.

Let us illustrate the two-step MLE-process. Now we take $N=n^{3/8}$.

\begin{figure}[h]
\centering\includegraphics [scale=0.68]  {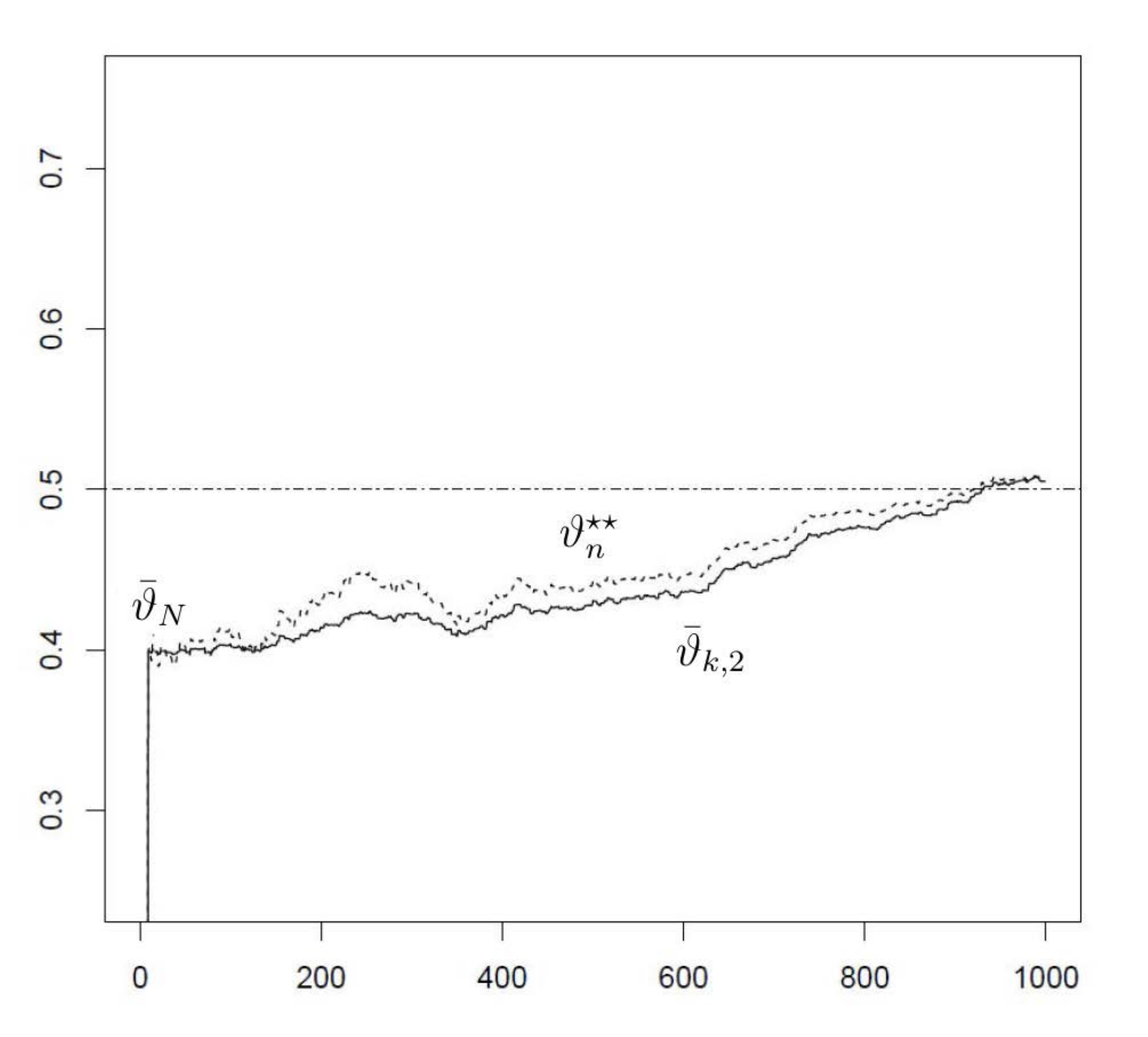}
\caption{Second preliminary and two-step MLE-processes. $n=10^3 $,  $\vartheta=0.5$}
\label{figure_22}
\end{figure}

\begin{figure}[h]
\centering\includegraphics [scale=0.68]  {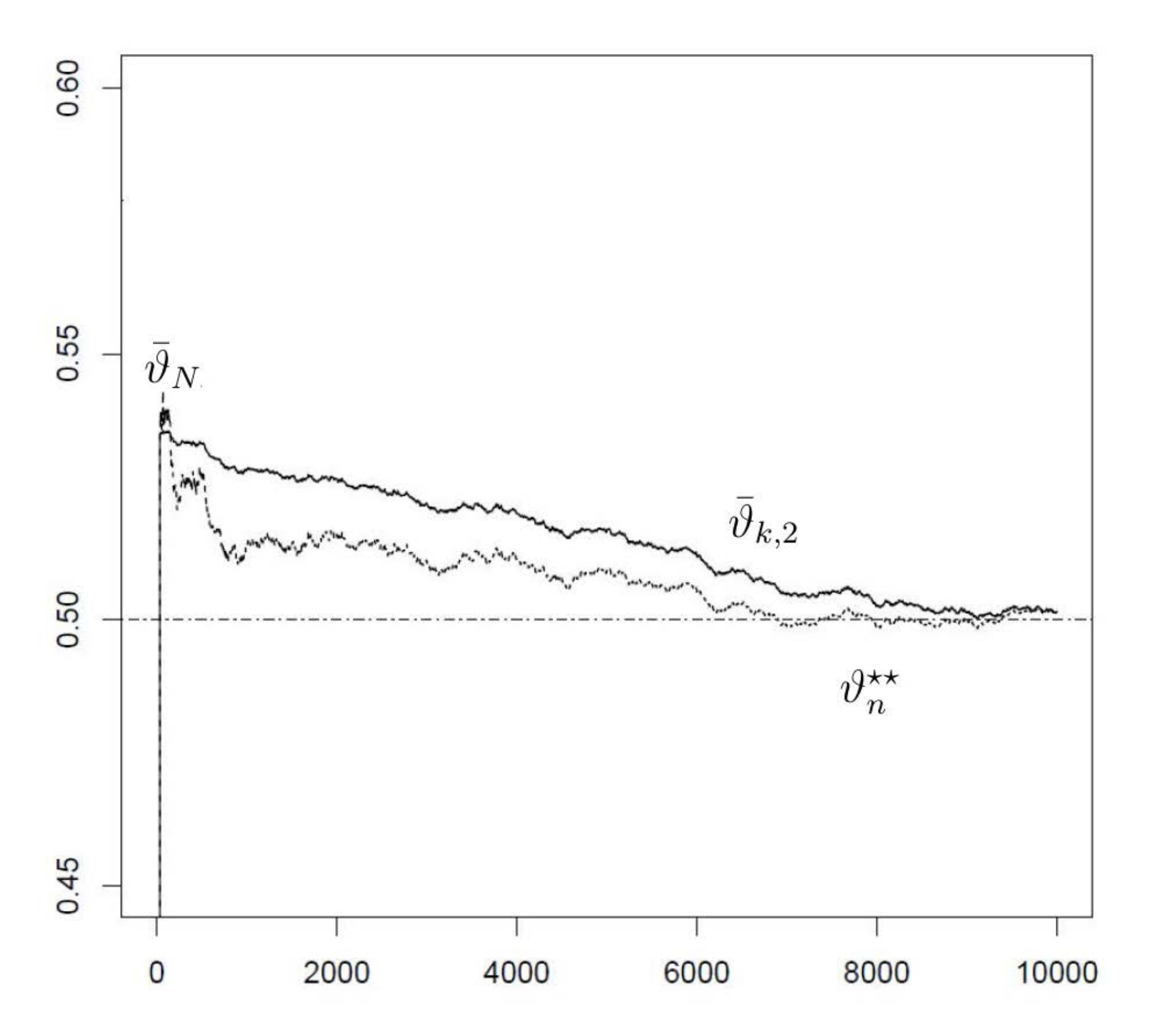}
\caption{Second preliminary and two-step MLE-processes. $n=10^4$, $\vartheta=0.5$}
\label{figure_23}
\end{figure}

We consider two cases: one with $n=10^3$ observations and the second  with
$n=10^4$ observations.

On the Figure \ref{figure_22} the preliminary EMM $\bar\vartheta _N=0.4$ that
is  far from the true value $\vartheta =0.5$. We obtain this estimator
based on the learning interval of $N=1000^{3/8}\approx 13$ observations. Then we obtain the
second preliminary estimator-process   $\bar\vartheta _{k,2}, 
k=N+1,\ldots,n)$ (continuous line) and see  that it  tends to the true value. The two-step
MLE-process $\vartheta _n^{\star \star}$ (dashed line) is closer to the true value and
as well tends to the true value.

On the Figure \ref{figure_23} the preliminary EMM $\bar\vartheta _N=0.54$ that
is close to the true value $\vartheta =0.5$. We obtain this estimator
based on the learning interval of $N=10000^{3/8}=32$ observations.
Then we obtain the
second preliminary estimator-process   $\bar\vartheta _{k,2}, 
k=N+1,\ldots,n)$ (continuous line) and see  that it  tends to the true value. The two-step MLE-process $\vartheta _n^{\star \star}$ (dashed line) is closer to the true value and as well tends to the true value.

\section{Discussion}

Two-step MLE-process allows us to estimate the parameter $\theta $ for the
values $k$ satisfying the condition $n^{1/4}<k\leq n$. If we need a shorter
  learning interval, say, $\left[1,n^\delta \right]$ with $\delta \in
  (\frac{1}{8}, \frac{1}{4}]$, then we have to study the three-step
    MLE-process, i.e., we use a preliminary estimator $\bar\vartheta _N$ and
    two estimator-processes like \eqref{kk}.

Note that the proposed one-step MLE-process can be written in the recurrent
form. Indeed, the estimator $\vartheta _{k,n}^\star$ we can write as follows
\begin{align*}
&\vartheta _{k+1,n}^\star=\bar\vartheta _N+\frac{1}{\sqrt{k+1}}\II\left(\bar\vartheta
_N\right)^{-1}\Delta _{k+1}\left(\bar\vartheta _N,X^{k+1} \right) \\
&\quad =\bar\vartheta _N+\frac{1}{{k+1}}\II\left(\bar\vartheta
_N\right)^{-1}\left[\sum_{j=1}^{k}\dot\ell\left(\bar\vartheta
_N, X_{j-1},X_j \right)+\dot\ell\left(\bar\vartheta
_N, X_{k},X_{k+1}\right)\right]\\
&\quad = \frac{k}{k+1}\left[ \bar\vartheta _N+\frac{1}{{k}}\II\left(\bar\vartheta
_N\right)^{-1}\sum_{j=1}^{k}\dot\ell\left(\bar\vartheta
_N, X_{j-1},X_j \right)\right]+\frac{1}{k+1} \bar\vartheta _N\\
&\qquad \quad +\frac{1}{{k+1}}\II\left(\bar\vartheta
_N\right)^{-1}\dot\ell\left(\bar\vartheta
_N, X_{k},X_{k+1}\right)\\
&\quad = \frac{k}{k+1}\vartheta _{k,n}^\star+\frac{1}{k+1} \bar\vartheta
  _N+\frac{1}{{k+1}}\II\left(\bar\vartheta 
_N\right)^{-1}\dot\ell\left(\bar\vartheta
_N, X_{k},X_{k+1}\right).
\end{align*}
The obtained presentation 
\begin{align*}
\vartheta _{k+1,n}^\star= \frac{k}{k+1}\vartheta _{k,n}^\star+\frac{1}{k+1}
\bar\vartheta _N+\frac{1}{{k+1}}\II\left(\bar\vartheta
_N\right)^{-1}\dot\ell\left(\bar\vartheta _N, X_{k},X_{k+1}\right)
\end{align*}
allows us to calculate $\vartheta _{k+1,n}^\star$ using the values
$\bar\vartheta _N,\vartheta _{k,n}^\star $ and observations $X_k,$ $X_{k+1}$ only. 

The similar structure can be obtained for the two-step MLE-process too.

{\bf Acknowledgment.} We are very grateful to the Referees for the comments which
allowed us to improve the exposition.  This work was done under partial financial support of
the grant of  RSF number 14-49-00079.


\begin{thebibliography}{99}

\bibitem {Bil61} Billingsley, P. (1961) {\it Statistical Inference for Markov
  Processes.} Chicago: The University of Chicago Press.

\bibitem {BSV15} Bogachev, V.I., Shaposhnikov, S.V. and Veretennikov,
  A.Yu. (2015) Differentiability of solution of stationary
  Fokker-Planck-Kolmogorov equations with respect to a parameter. {\it Doklady
    Mathematics}, 91, 1, 76-79.

\bibitem{FH82} Fabian, V., Hannan, J. (1982) On estimation and adaptive
  estimation for LAN families. {\it Z. Wahrsch. Verw. Geb.}, 59, 459-479.


\bibitem{FY03} Fan, J. and Yao, Q. (2003) {\it Nonlinear Time Series:
  Nonparametric and Parametric Methods}. Springer-Verlag, New York.

\bibitem {IH81} Ibragimov I.A. and Khasminskii R. (1981) {\it Statistical
Estimation. Asymptotic Theory.} {Springer-Verlag}, New York.

\bibitem{KU15} Kamatani, K. and Uchida, M. (2015) { Hybrid multi-step
  estimators for stochastic differential equations based on sampled data.}
  {\it Statist. Inference Stoch. Processes.} 18, 2, 177-204.  

	\bibitem {K14}  Kutoyants, Y.A. (2014) {On approximation of the
  backward stochastic differential equation. Small noise, large samples and
  high frequency cases.}  {\it Proc. Steklov Inst. Mathem.}, 287, 133-154.
\bibitem {Kut15} Kutoyants, Y.A. (2015) On multi-step MLE-processes for
  ergodic diffusion. submitted.
\bibitem {KZ1} Kutoyants, Y.A. and Zhou, L. (2014) {On approximation of the
  backward stochastic differential equation. }
   {\it J. Stat. Plann. Infer.  150, 111-123}.
\bibitem{LC56} Le Cam, L. (1956) On the asymptotic theory of estimation and
  testing hypotheses. {\it Proc. 3rd Berkeley Symposium I}, 355-368.

\bibitem {MA15}     Motrunich, A. (2015) On parameter estimation for Markov
  sequences. Submitted. 
\bibitem{OI77} Ogata, Y. and Inagaki, N. (1977) The weak convergence of the
  likelihood ratio random fields for Markov observations. {\it
    Ann. Inst. Statist. Math.}, 29, Part A, 165-187.

\bibitem{RG65} Roussas, G.G. (1965)  Asymptotic inference in Markov
  processes. {\it Ann. Math. Statistics}, 36, 3, 978-992.

\bibitem{RG72} Roussas, G.G. (1972) {\it Contiguity of Probability Measures: Some
  Applications in Statistics}. Cambridge Univ. Press, London.

\bibitem{SKh96} Skorohod, A.V. and Khasminskii, R.Z. (1996) On parameter
  estimation by indirect observations. {\it Prob. Inform. Transm.}, 32, 58-68.

\bibitem{TK00} Taniguchi , M. and Kakizawa, Y. (2000) {\it Asymptotic Theory
  of Statistical Inference for Time Series.} {Springer-Verlag}, New York.

\bibitem {VV02} Varakin A.B., Veretenikov A.Yu. (2002) {\sl On parameter
  estimation for "polynomial ergodic" Markov chains with polynomial growth
  loss functions.} Markov Processes and Related Fields, 8(1), 127-144.

\bibitem{Ver00} Veretenikov, A.Yu. (2000) {\it Parametric and Non-Parametric
  Estimation of Markov Chains} [in Russian]. Moscow Univ. Press, Moscow.

\end{thebibliography}
\end{document}